\documentclass[twoside]{article}
\usepackage{stmaryrd}

\usepackage[all]{xy}
\usepackage{amsfonts}
\usepackage{amsmath}
\usepackage{amssymb}
\usepackage{latexsym}
\usepackage{mathrsfs}

\newtheorem{thm}{\bf Theorem}[section]
\newtheorem{cor}[thm]{\bf Corollary}
\newtheorem{lem}[thm]{\bf Lemma}
\newtheorem{prop}[thm]{\bf Proposition}

\def \pde{\rm pd_{{\rm\underline{\underline{A}}}\ltimes_n F}}
\def \pda{\rm pd_{{\rm\underline{\underline{A}}}}}
\def \ide{\rm  id_{G\ltimes_n {\rm\underline{\underline{A}}}}}
\def \ida{ \rm id_{{\rm\underline{\underline{A}}}}}
\def \idrn{\rm id_{R\ltimes_{n}M}  }
\def \idr{\rm id_{R}}

\def\pr{{\parindent0pt {\bf Proof.\ }}}

\def\Fij{{\rm F_{i+j}}}
\def\Fjk{{\rm F_{j+k}}}
\def\Fi{{\rm F_i}}
\def\Fj{{\rm F_j}}
\def\Fk{{\rm F_k}}

\def\FI{{\rm (F_i)_{i\in I}}}
\def\Gij{{\rm G_{i+j}}}
\def\Gj{{\rm G_j}}
\def\Gi{{\rm G_i}}
\def\Gk{{\rm G_k}}
\def\Gjk{{\rm G_{j+k}}}
\def\Gijk{{\rm G_{i+j+k}}}
\def\GI{{\rm (G_i)_{i\in I}}}
\def\ext{{\rm \underline{\underline{A}}\ltimes F}}
\def\A{\rm \underline{\underline{A}}}
\def\Rn{ R\ltimes_{n} M }
\def\Rm{ R\ltimes M }

\def\RM{ R\ltimes_{n} M_{1}\ltimes \cdots\ltimes M_{n} }
\def\lRM{ _{R}Mod\ltimes_{n} (M_{i}\otimes -)_i }
\def\RRM{  (\Hom(M_{i},-))_i\rtimes_{n} Mod_{R} }

\def \ntex{ {\rm\underline{\underline{A}}\ltimes_{n} F } }

\def\conext{{\rm G\rtimes _{n} \underline{\underline{A}}    }     }
\def\tex{ {\rm\underline{\underline{A}}\ltimes_{\Phi}F_{1}\ltimes\cdots\ltimes F_{n}}}
\def\cotex{{\rm G_{n}\rtimes\cdots\rtimes G_{1}\rtimes_\Psi \underline{\underline{A}}} }
\def\Im{{\rm Im}}
\def\Coker{{\rm Coker}}
\def\Ker{{\rm Ker}}
\def\Ext{{\rm Ext}}

\def\Hom{{\rm Hom}}

\newcommand{\cqfd}
{\hspace{1cm}
\rule{2mm}{2mm}%
\medbreak%
\par%
}

\begin{document}

\title{The category of modules on an $n$-trivial extension: the basic properties
}

\author{Dirar Benkhadra$^{1,a}$,  Driss Bennis$^{1,b}$  and J. R. Garc\'{\i}a Rozas$^2$}

\date{}

\maketitle

\small{1: CeReMAR Center,   Faculty of Sciences,    Mohammed V University in Rabat,  Morocco.\\
\noindent $^a$ dirar.benkhadra@hotmail.fr; $^b$ driss.bennis@um5.ac.ma; driss$\_$bennis@hotmail.com

2: Departamento de  Matem\'{a}ticas, Universidad de Almer\'{i}a, 04071 Almer\'{i}a, Spain.

\noindent e-mail address: jrgrozas@ual.es
 }

\date{}

\maketitle

\begin{abstract}
In this paper we  investigate a categorical aspect of $n$-trivial extension of a ring by a family of modules. Namely, we introduce the right (resp.,
left) $n$-trivial extension of a category by a family of endofunctors.   Among other results, projective, injective and flat objects of this category
are characterized. We end the paper with two applications. We  characterize when  an $n$-trivial extension ring is $k$-perfect   and we establish a
result on the selfinjective dimension of an $n$-trivial extension ring.

 \end{abstract}



\section{Introduction}
Throughout this paper $R$ will be an associative (non necessarily commutative) ring with identity, and all modules will be, unless otherwise specified, unital left $R$-modules.  The category of all left (right if needed) $R$-modules will be denoted by $_{R}Mod$  ($Mod_{R}$).

The notion of trivial extension of a ring by a bimodule has been played an important role in various parts of algebra. For instance, in ring theory after the work of Nagata \cite{Nag} in which this construction is called  ``idealization". In \cite{TEx}, Fossum, Griffith and Reiten investigated  the categorical aspect of the trivial extension and gave a detailed treatment of its homological properties.\\

Recently, in \cite{Ben}, the classical trivial extension has been extended by associating a ring to a ring $R$ and a family $M=(M_i)^{n}_{i=1} $ of $n$ $R$-bimodules for a given $n\geqslant 1$, which is called $n$-trivial extension of $R$ by $M$ and denoted by $\Rn$.  Precisely,  for a ring $R$ and a family $M=(M_i)^n_{i=1}$  of $R$-bimodules and $\Phi=(\Phi_{i,j})_{1\leqslant i,j\leqslant n-1}$ a family of bilinear maps such that each $\Phi_{i,j}$ is written multiplicatively $$\Phi_{i,j}: \underset{(m_i,m_j)}{M_i\times M_j} \longrightarrow \underset{\Phi(m_i,m_j):=m_i.m_j}{M_{i+j}}$$ whenever $i+j\leq n$  with, $$\Phi_{i+j,k}(\Phi_{i,j}(m_i,m_j), m_k)=\Phi_{i,j+k}((m_i,\Phi_{i,k}(m_j,m_k)). $$

\par The $n$-trivial extension ring of $R$ by $M=(M_i)^n_{i=1}$ is the additive group $R\oplus M_1\oplus \cdots \oplus M_n$ associated with the multiplication: $(m_0,...,m_n)(m'_0,...,m'_n)=\sum_{i+j=n}m_im'_j$. Notice that, when $n=1$, $\Rn$ is nothing but the classical trivial extension of $R$ by $M$ which is denoted by  $\Rm$. During this paper, $M=(M_i)^n_{i=1}$ is a family of $R$-bimodules with the above assumptions on $\Phi_{i,j}\,'s$, so   when using  $\Rn$ we mean the $n$-trivial extension of $R$ by  $M=(M_i)^n_{i=1}$. \\

In   \cite{Ben}, various ring-theory proprieties of  $n$-trivial extensions were studied.  In this paper, we are going to study some categorical
aspects of the category of modules over an $n$-trivial extension. Our work is based on the work  of Fossum, Griffith and Reiten \cite{TEx} on the so
called  right  trivial extension of a  category by an endofunctor. Recall from \cite{TEx} that the classical right trivial extension of a category
$\A$ by an endofunctor $F$  is a category whose objects are  couples $(X,f)$ with $X\in ob(\A)$ and $f:FX\longrightarrow X$ such that $Ff\circ f=0$;
and a morphism
    $\gamma :(X,\alpha) \longrightarrow (Y,\beta) $     is a      morphism $X \longrightarrow Y$ in $\A$  such that the diagram
\begin{displaymath}
\xymatrix {F X\ar[d]^{\alpha}\ar[r]^{F\gamma}& F Y\ar[d]^{\beta}\\
X \ar[r]_{\gamma}&Y }
\end{displaymath} commutes. This category is denoted by $\ext$. It is proved that, when $\A$ is the category of left $R$-modules and $F$ is the tensor product
functor $M\otimes_{R}-$, $\ext$ is isomorphic to the category of left modules over $\Rm$.\\

Here we introduce a general concept which we call the right $n$-trivial extension of $\A$ by a family of functors $F=(F_i)^n_{i=1}$.  This is done,
following the way adopted in \cite{TEx} (see also \cite{Cl}),  by carefully analyzing the left module action on an $n$-trivial extension.  Then, we investigate some of its   properties including the homological ones. \\

We have structured the paper in the following way:\\

In Section 2, we start with a paragraph in which we define the right $n$-trivial extension of a  category $\A $ by a family of endofunctors $F=(\Fi)^n_{i=1}$, $\ntex$.
 When, $n=1$, $\ntex$ is nothing but the classical trivial extension defined in \cite{TEx}.  Then, we show when the category $\ntex$ is abelian (see Proposition  \ref{prop.ab1}).  For  the case where $\A$ is the category of  left  $R$-modules and  $F_i=M_i\otimes_R-$ for each $i\in \{1,...,n\}$, where $M_i$ are  $R$-bimodules, $\ntex$ is nothing but the category of left modules over $\Rn$ (Theorem \ref{thm.trivmod}).\\
Furthermore, we define the left   $n$-trivial extension of a  category $\A $ by a family of endofunctors $G=(\Gi)^n_{i=1}$, $\conext$.
We  show   when it is an abelian category (Proposition  \ref{prop.ab2}) and when is isomorphic to a right   $n$-trivial extension (Theorem \ref{thm.RtrivLtriv}).\\
%%%%%%%%%%%%%%%%%%%%%%%%%%%%%%%%%%%%%%%%
%%%%%%%%%%%%%%%%%%%%%%%%%%%%%%%%%%%%%%%%%

Section 3 is devoted to the study of some homological notions. We begin by characterizing the class of objects which are isomorphic to direct summands of direct sums of copies of an object   in $\ntex$ (see Theorem \ref{thm.Add}). This leads to a  characterisation of  projective  objects in $\ntex$ (see Corollary \ref{cor.proj}). Dually, a     characterisation of  injective is given (see Theorem \ref{thm.Prod} and Corollary \ref{cor.inj}). Flat objects are also studied  (see Theorem \ref{thm.Lim}, Corollaries \ref{cor.flat} and \ref{cor.flatmod}).\\

 Finally in Section 4  we establish   some results as applications of the previous study. Namely, we  characterize when  the $n$-trivial extension $\Rn$ is $k$-perfect (Theorem \ref{thm.perfect}) and we establish a result on the selfinjective dimension of the ring $\Rn$ (Theorem \ref{thm.selfinj}).

\bigskip

%%%%%%%%%%%%%%%%%%%%%%%%%%%%%%%%%%%%%%%%%%%%%%%%%%%%%%%%%
%%%%%%%%%%%%%%%%%%%%%%%%%%%%%%%%%%%%%%%%%%%%%%%%%%%%%%%%%
%%% %%%%%%%%%%%%%%%%%%%%%%%%%%%%%%%%%%%%%%%%%%%%%

%%%%%%%%%%%%%%%%%%%%%%%%%%%%%%%%%%%%%%%%%%%%%%%%%%%%%%%%%
%%%%%%%%%%%%%%%%%%%%%%%%%%%%%%%%%%%%%%%%%%%%%%%%%%%%%%%%%
%%%%%%%%%%% %%%%%%%%%%%%%%%%%%%%%%%%%%%%%%%%%%%%%%%%%%%%
\section{n-Trivial extensions of abelian categories}
As mentioned in the introduction, we aim to establish some basic properties of the category of modules over $S=\RM$. To this end we follow the way
adopted in \cite{TEx} and \cite{Cl}. Explicitly, we introduce  a category which is equivalent to  the category of left modules $\RM$-mod. The definition of the desired category
follows from the analysis of the module action on an $n$-trivial extension. To this end, we need to consider the tensor products rather than   pre-products $\Phi_{i,j}$. In fact, by \cite{Ben}, any pre-product maps $\Phi_{i,j}: M_i\times_{R}M_j\longrightarrow M_{i+j}  $ have a corresponding $R$-module homomorphism
$$\tilde{\Phi}_{i,j}: M_i\otimes_{R}M_j\longrightarrow M_{i+j}  $$ with $i+j\leq n$ i.e $$\tilde{\Phi}_{i+j,k}\circ(\tilde{\Phi}_{i,j}\otimes_R id_{M_k})=\tilde{\Phi}_{i,j+k}\circ(id_{M_i}\otimes_R \tilde{\Phi}_{i,k}).$$ We set $m_im_j:=\Phi_{i,j}(m_i, m_j)=\tilde{\Phi}_{i,j}(m_i\otimes m_j)$.\\
Now we consider an $S$-module $X$. Via the natural injection to the first coordinate,  $i:R\longrightarrow S$, $X$ is an $R$-module.
Thus $S$ as a right $R$-module gives rise to a left homomorphism of $R$-modules
 \begin{displaymath}
    \xymatrix{ S\otimes_R X\ar[r]^{u}& X:
                     s\otimes x\ar@{|->}[r]& s.x     }
\end{displaymath}
which satisfies the commutativity of the two following diagrams:

 \begin{displaymath}
    \xymatrix{
           R\otimes_R X \ar[r]^{i\otimes id_{X}} \ar[rd]_{=} & S\otimes_{R} X \ar[d]^{u} \\
              & X }
\end{displaymath}
 and \begin{displaymath}
    \xymatrix{
        S\otimes_R S\otimes_R X \ar[r]^{m\otimes X} \ar[d]_{S\otimes_R u} &S\otimes_R X\ar[d]^{u} \\
        S\otimes_R X\ar[r]_{u}& X }
\end{displaymath} \\

It is also clear that the above data determines uniquely the $S$-module structure of $X$.\\
As an $R$-module $S\otimes_R X\cong X\oplus (M_1\otimes X)\oplus \cdots\oplus (M_n\otimes X)$, so $u=(r,\alpha_1,\dots,\alpha_n)$, where $\alpha_i:M_i\otimes X\longrightarrow X$ is a homomorphism of $R$-modules for all $i\in I$ and $r$ is an endomorphism of $X$. By the first diagram $r=id_X$ and by the second one we have: \begin{itemize}
     \item[i)] $\alpha_i \circ(M_i\otimes\alpha_j)=0$ for $i+j>n.$\\

    \item[ii)]  The following diagram \begin{displaymath}
    \xymatrix{
       M_{i}\otimes M_{j}\otimes X \ar[r]^{ \quad \Phi_{i,j}\otimes id{_X}} \ar[d]_{ id_{M_{1}}\otimes \alpha_{j}} & M_{i+j}\otimes X  \ar[d]^{\alpha_{i+j}} \\
      M_{i}\otimes X \ar[r]_{\alpha_{i}} & X      }
\end{displaymath} commutes.
 \end{itemize}
On the other hand, if we take a left $R$-module $X$ and a family $(\alpha_i:M_i\otimes X\longrightarrow X)_{i\in I}$ verifying the two conditions
$(i)$ and $(ii)$ above, we could define a left structure of $S$-module via the module  action:
$$(a,m_1,m_2,...,m_n).x=ax+\Sigma_{i=1}^{n}\alpha_i(m_i\otimes x ).$$

This leads as done in \cite{TEx} to introduce the category of right $n$-trivial extension of an abelian category  as follows:\\
Let $\A$ be an abelian category, $I=\{ 1,2,\cdots , n\}$ an index set, $\rm{F}=\FI$ be a family  of additive covariant endofunctors of $\A$, and $
\Phi=(\Phi_{i,j} : \Fi \Fj \longrightarrow \Fij )$  be a family of natural transformations such that, for every $i,j,k\in I$ with $i+j+k\in I$, the
diagram:
\begin{displaymath}
\xymatrix {\Fi \Fj \Fk X \ar[r]^{\Fi (\Phi_{j,k})_X}\ar[d]_{(\Phi_{i,j })_{\Fk X}}&  \Fi \Fjk X\ar[d]^{(\Phi_{i,j})_X}\\
              \Fij\Fk X \ar[r]_{(\Phi_{i+j,k })_X}& F_{i+j+k} X}
\end{displaymath} commutes.\\
 We define an additive category which we denote by  $\tex$ (or simply by $\ntex$) and we call it the right $n$-trivial extension of $\A$ by $\rm{F}$, as follows:

    \textbf{Objects}:  An object of  $\tex$ is of the form  $(X,f)$ with $X$  an object of $\A$ and $f=(f_i)_{i\in I}$  a family of morphisms $f_i:\Fi X\longrightarrow X$ such that for every $i,j\in I$, the diagram
    \begin{displaymath}
\xymatrix {\Fi \Fj X\ar[d]_{\Fi f_j}\ar[r]^{(\Phi_{i,jX})_X}& \Fij X\ar[d]^{f_{i+j}}\\
                \Fi X \ar[r]_{f_i}& X}
\end{displaymath} commutes when $i+j\in I$, if not, the composition
 \begin{displaymath}
\xymatrix {\Fi \Fj X \ar[r]^{\Fi f_j}&\Fi X\ar[r]^{f_i}& X}
\end{displaymath} is zero.

    \textbf{Morphisms}: For every two objects $(X,\alpha)$ and $(Y,\beta)$, a morphism
 $\gamma :(X,\alpha) \longrightarrow (Y,\beta) $ in $\ntex$  is   a morphism $X \longrightarrow Y$ in $\A$  such that, for every $i\in I$, the diagram
\begin{displaymath}
\xymatrix {\Fi X\ar[d]_{\alpha_i}\ar[r]^{\Fi\gamma}& \Fi Y\ar[d]^{\beta_i}\\
X \ar[r]_{\gamma}&Y }
\end{displaymath} commutes.

One can remark that, for $n=1$, \underline{\underline{A}}$\ltimes_1$F is nothing but the classical trivial extension of \cite{TEx}.\\
It is clear that $\ntex$ is an additive category. Following the same argument as done in \cite[Proposition 1.1]{TEx}, we can show the following
result.

\begin{prop}\label{prop.ab1}
If $\Fi$ is right exact for every $i\in I$, then $\ntex$   is an abelian category.
\end{prop}
One can notice that the description of modules over  $n$-trivial extension $\Rn$  given in the first part of this section leads to our first main purpose in this section.
\begin{thm}\label{thm.trivmod}
The category  $\lRM$ is isomorphic to the category of left modules over $\RM$.
\end{thm}
\par Furthermore the right $n$-trivial extension can be seen as a  cleft extension which is introduced by Beligiannis in \cite{Cl} as follows:
 A category $\mathcal{D}$ is called a right cleft extension of a category $\mathcal{A}$ if there exists a tripleable functor $U:\mathcal{D}\longrightarrow \mathcal{A}$ and a functor $Z:\mathcal{A}\longrightarrow \mathcal{D}$ such that $UZ\cong id_{\mathcal{A}}$.\\
In order to show that  $\ntex$ is a right cleft extension, we define $(U,Z)$ as follows:\\
The underlying functor $U: \ntex  \longrightarrow  \A$ is defined by $U(X,f)=X$ for any object $(X,f)$ in $\ntex$   and for an arrow $\alpha$ in $\ntex$, $U(\alpha)=\alpha$.\\
The zero functor $Z:\A\longrightarrow\ntex$ is defined by $Z(X):(X,0:\Fi X\shortrightarrow X)$ and for morphisms $\alpha:X\longrightarrow Y$ in $\A$, $Z(\alpha)=Z(X)\longrightarrow Z(Y)$.\\

 Also, as in \cite{Cl}, we can construct   two other functors $T$ and $C$ as follows:\\
 The functor $T: \A  \longrightarrow \ntex$  given, for every object $X$, by $T(X)=(X\oplus(\displaystyle\bigoplus_{i\in I}\Fi X),\kappa)$, where $\kappa=(\kappa_i)$ such that,
for every $i\in I$,   $$\kappa_i=(a^i_{\alpha,\beta}):\Fi X\displaystyle\bigoplus_{j\in I}\Fi\Fj X\longrightarrow X\oplus(\displaystyle\bigoplus_{j\in I}\Fj X)$$
defined by\\
\begin{equation*}
  \begin{cases}
 a^i_{i+1,1}=id_{\Fi X}.\\
a^i_{i+j+1,j+1}=\Phi_{ij X} & \text{for j with}\quad  i+j\in I.\\
    0 &\text{otherwise}.
  \end{cases}
\end{equation*}
 and for morphisms by:
$$T(\alpha)=\left( \begin{array}{cccc}
 \alpha &0  &\ldots &0\\
 \vdots &F_1\alpha &   & \vdots\\
 0   & \ldots &   \ddots&0\\
 0 & \ldots&  0& F_{n}\alpha
\end{array}\right)$$
 The cokernel functor  $C:\ntex\longrightarrow \A$  is given by  $C((X,f))=\Coker(f:\oplus_{i\in I}\Fi X\longrightarrow X)$  for every object $(X,f)$ and for a morphism $\alpha:(X,f)\longrightarrow (Y,g)$ in $\ntex$, $C(\alpha)$ is the induced map.\\

We can verify that $(T,U)$ and $(C,Z)$ are pairs of adjoint functors with $CT\cong id_{\A}$ and $UZ\cong id_{\A}$. These functors have a concrete realization on $\Rn$.\\

 One can notice that, for every $\Rn$-module $X$, $U(X)$ is nothing but the  $R$-module $X$ (through $i:R\rightarrow \Rn$), while $CX=X\slash  \Sigma_{i\in I} M_i X$. On the other hand, for an $R$-module X,  $TX=\Rn\otimes_{\A} X$ and  $ZX$ is the  $\Rn$-module defined through the projection $\pi: \Rn\longrightarrow R$. \\

%%%%%%%%%%%%%%%%%%%%%%%%%%%%%%%%%%%%%%%%%%%
%%%%%%%%%%%%%%%%%%%%%%%%%%%%%%%%%%%%%%%%%%%%%%%%%%%%%
%%%%%%%%%%%%%%%%%%%%%%%%%%%%%%%%%%%%%%%%%%%%%%%%%%%%%%%%%%
%%%%%%%%%%%%%%%%%%%%%%%%%%%%%%%%%%%%%%%%%%%%%%%%%%%%%%%%%%%%%%%
%%%%%%%%%%%%%%%%%%%%%%%%%%%%%%%%%%%%%%%%%%%%%%%%%%%%%%%%%%%%%%%%%
\par To give the ``categorical" aspect of the category of right modules over an  $n$-trivial extension, we use as done in \cite{TEx} the functors $\Hom_R(M_i,-)$. This leads to the generalisation of the classical left trivial extension of abelian categories.

\par Let G$=\GI$ be a family of additive covariant endofunctors of $\A$
and   $\Psi=(\Psi_{i,j}:\Gij X\longrightarrow \Gj\Gi)$ be a family of natural transformations such that, for every $i,j,k\in I$ with $i+j+k\in I$,
\begin{displaymath}
\xymatrix { \Gi\Gjk X \ar[r]^{\Gi(\Psi_{j,k})_X} & \Gi\Gj\Gk X\\
               \Gijk X\ar[u]^{(\Psi_{i,j+k})_X}\ar[r]_{(\Psi_{i+j,k})_X}& \Gij\Gk X\ar[u]_{(\Psi_{i,j})_{\Gk X}} }
\end{displaymath} commutes.\\
\\
 We define the left $n$-trivial extension of $\A$   by G denoted by $\cotex$ or simply $\conext$ as follows:\\
 \textbf{Objects}: An object of   $\conext$  is a couple $(X,g)$, where $X$ is in $\A$ and $g=(g_i)_{i\in I}$ is a family of morphisms
 $g_i:X \longrightarrow \Gi X$ such that for $i\in I$ the diagram
\begin{displaymath}
\xymatrix { X \ar[d]_{g_{i+j}}\ar[r]^{g_i} & \Gi X\ar[d]^{\Gi g_j}\\
                \Gij X \ar[r]_{\Psi_{ij}}&\Gi \Gj X}
\end{displaymath} commutes when $i+j\in I$, if not, the composition
 \begin{displaymath}
\xymatrix { X\ar[r]^{g_i} & \Gi X \ar[r]^{ \Gi g_j}&\Gi\Gj X \\}
\end{displaymath} is zero.\\
\textbf{Morphisms}: For every two objects $(X,\alpha)$ and $(Y,\beta)$ in $\cotex$, a morphism
  $\gamma :(X,\alpha) \longrightarrow (Y,\beta)$ is   a morphism  $X \longrightarrow Y$ in   $\A$ such that the diagram
\begin{displaymath}
\xymatrix {\Gi X \ar[r]^{\Gi\gamma}& \Gi Y\\
X \ar[r]_{\gamma} \ar[u]^{\alpha_i}&Y \ar[u]_{\beta_i}}
\end{displaymath} commutes for every  $i\in I$.

Similarly to the right $n$-trivial extension we have the following results.

\begin{prop}\label{prop.ab2}
The category $\conext$ is abelian if $\Gi$ is left exact for every $i\in I$.
\end{prop}

\begin{thm}\label{thm.rtrivmod}
The category  $\RRM$ is isomorphic to the category of right modules over $\RM$.
\end{thm}

In the context of left $n$-trivial extensions, we have pairs of adjoint functors $(U,H)$ and $(Z,K)$ defined as follows:

The functor  $H:\A \longrightarrow\conext$ given, for every $X\in \A$, by $H(X)=(G_nX\oplus G_{n-1}X\oplus\cdots \oplus G_1 X)\oplus X,\lambda=(\lambda_i))$, where,  for every $i\in I$,
$$\lambda_{i}=(a^i_{\alpha,\beta}):G_iG_nX\oplus G_i G_{n-1}X\oplus\cdots \oplus G_iG_1 X\oplus G_iX\longrightarrow G_nX \oplus G_{n-1}X\oplus\cdots \oplus G_1 X \oplus X $$ defined by
\begin{equation*}
  \begin{cases}
 a^i_{n+1,n+1-i}=id_{\Gi X}.\\
a^i_{n-j+1,n+1-(i+j)}=\Psi_{ijX} & \text{for j where}\;  i+j\in I.\\
    0 &\text{otherwise}.
  \end{cases}
\end{equation*}
and on morphisms by $$H(\alpha)=\left( \begin{array}{cccc}
 G_{n}\alpha &0  &\ldots &0\\
 \vdots &\ddots &   & \vdots\\
 0   & \ldots &G_1\alpha   &0\\
 0 & \ldots&  0& \alpha
\end{array}\right)$$
 The underlying functor    $U:\conext  \longrightarrow \A$ is defined by $U(X,g)=X$ for any object in $\conext$  and for an arrow $\alpha$ in $\conext$, $U(\alpha)=\alpha$.
 The zero functor $Z:\A\longrightarrow\conext$ is defined by $Z(X)=(X,0: X\shortrightarrow \Gi X)$ and for morphism $\alpha:X\longrightarrow Y$ in $\A$ we get  $Z(\alpha)=Z(X)\longrightarrow Z(Y) $ in $\conext$.
 Finally  the kernel functor $K:\conext\longrightarrow\A$  is given, for every object in $(\cotex)$ by $K(X,g:X\longrightarrow \Gi X)=\Ker(g:X\longrightarrow\oplus_{i\in I} \Gi X )$, and for a morphism $\beta$, $K\beta$ is the induced morphism.\\

We easily check that $UZ\cong id_{\A}$ and $KH\cong id_{\A}$.\\

\par In the concret case of the category of left modules over the  $n$-trivial extension ring, the functor $H$ is simply the ``hom" functor $Hom_{R}(S,-)$ with $S=\RM$.\\

%%%%%%%%%%%%%%%%%%%%%%%%%%%%%%%%%%%%%%%%%%%%%%%%%%%%%%%%%
%%%%%%%%%%%%%%%%%%%%%%%%%%%%%%%%%%%%%%%%%%%%%%%%%%%%%%%%%
%%%%%%%%%%%%%%%%%%%%%%%%%%%%%%%%%%%%%%%%%%%%%%%%%%%%%%%
 \par We end this section with a generalization of \cite[Proposition 1.11]{TEx}  which investigates the relation between $\conext$  and $\ntex$.\\

\par  In what follows we suppose that $\Fi$ admits a right adjoint functor $\Gi$ for every $i\in I$, so we have a natural isomorphism
\begin{displaymath}
    \xymatrix {\Hom(\Fi X,Y)\ar[r]^{\phi_i}&\Hom(X,\Gi Y). }
\end{displaymath}
We recall from \cite[page 81]{Mac} that $$\begin{array}{cc}
\phi_i(f\circ \Fi h)=\phi_i(f)\circ h & \phi_{i}(k\circ f)=\Fi k\circ\phi_{i}(f).\\
\phi^{-1}_{j}(g\circ h)=\phi^{-1}_{j}(g)\circ \Gj h & \phi^{-1}_{j}(\Gj k\circ g)=k\circ\phi^{-1}_{j}(g)
\end{array} $$for $i,j\in I$,
where $f:\Fi X\shortrightarrow Y,$  $h:Z\shortrightarrow X,\quad k:Y \shortrightarrow W $ and $g:X\shortrightarrow \Gj Y$.

\begin{thm}\label{thm.RtrivLtriv}
Every right $n$-trivial extension of $\A$ by $F$, $\ntex$ is isomorphic to a left $n$-trivial extension of $\A$ by $G$, $\conext$ for  certain natural transformations.
\end{thm}

\pr For every object $(X,f)$ of $\ntex$, we consider the corresponding object of $\conext$ to be $(X,g)$ with $g=(\phi_i(f_i))_{i\in I}$ and an
associated natural transformation $$\Psi_X=\Gj\Gi\phi_{i+j}^{-1}(id_{\Gij X})\circ \phi_j\phi_i(\Phi_{(\Gij X)})$$ (here, since there is no
ambiguity, $\Phi_{ij}$ is   denoted simply by $\Phi$). We prove   that $(X,g)$ is an object of $\conext$. It is easy to show that $\Psi$ is a natural
transformation, so it remains to prove the commutativity of the diagram \begin{displaymath}
\xymatrix { \Gi\Gjk X \ar[r]^{\Gi\Psi_X} & \Gi\Gj\Gk X\\
               \Gijk X\ar[u]^{\Psi_X}\ar[r]_{\Psi_X}& \Gij\Gk X\ar[u]_{\Psi_{\Gk X}}         }
\end{displaymath}
We have:
\begin{eqnarray*}
    \Psi_X\circ g_{i+j}&=&\Gj\Gi\phi_{i+j}^{-1}(id_{\Gij X})\circ \phi_j\phi_i(\Phi_{ij}(\Gij X))\circ g_{i+j}\\
       &=&\Gj\Gi\phi_{i+j}^{-1}(id_{\Gij X})\circ \phi_j(\phi_i(\Phi_{\Gij X})\circ \Fj g_{i+j})\\
&=&\Gj\Gi\phi_{i+j}^{-1}(id_{\Gij X})\circ \phi_j\circ\phi_i(\Phi_{\Gij X}\circ \Fi\Fj g_{i+j})\\
&=&\phi_j[\Gi\phi_{i+j}^{-1}(id_{\Gij X})\circ \phi_i(\Phi_{\Gij X}\circ \Fi\Fj g_{i+j})]\\
&=&\phi_j\circ\phi_i[\phi_{i+j}^{-1}(id_{\Gij X})\circ \Phi_{\Gij X}\circ \Fi\Fj g_{i+j}]
\end{eqnarray*}
 Using the following diagram
 \begin{displaymath}
\xymatrix { \Fi\Fj X \ar[r]^{\Phi_X}\ar[d]_{\Fi\Fj g_{i+j}} & \Fij X\ar[d]^{\Fij g_{i+j}}\\
             \Fi\Fj\Gij X\ar[r]_{\Phi_{\Gij X}}& \Fij\Gij X }
\end{displaymath}
we get
\begin{eqnarray*}
\Psi_X\circ g_{i+j}&=&\phi_j\circ\phi_i[\phi_{i+j}^{-1}(g_{i+j})\circ \Phi_X]\\
&=& \phi_j\circ\phi_i[f_{i+j}\circ \Phi_X ]\\
&=& \phi_j\circ\phi_i(f_i\circ \Fi f_j)\\
&=& \phi_j(\phi_i(f_i)\circ f_j)\\
\Psi_X\circ g_{i+j}&=&\Gj g_{i}\circ g_j.
\end{eqnarray*}
Hence $(X,g)$ is an object of $\conext$.\\
Now, for an arrow $(X,f)\longrightarrow(X',f')$ in $\ntex$, we associate the induced arrow $(X,g)\longrightarrow(X',g')$ in $\conext$.\\
So far,  we have defined a functor  $\Delta: \ntex \longrightarrow \conext$ which is an isomorphism following standard arguments. \cqfd

%%%%%%%%%%%%%%%%%%%%%%%%%%%%%%%%%%%%%%%%%%%%%%%%
%%%%%%%%%%%%%%%%%%%%%%%%%%%%%%%%%%%%%%%%%%%%%%%%%%%%%%%%%%%%%%%%%%%%%%
%%%%%%%%%%%%%%%%%%%%%%%%%%%%%%%%%%%%%%%%%%%%%%%%%%%%%%%%%%%%%%%%%%%%%%%%%%%
%%%%%%%%%%%%%%%%%%%%%%%%%%%%%%%%%%%%%%%%%%%%%%%%%%%%%%%%%%%%%%%%%%%%%%%%%%%%%
%%%%%%%%%%%%%%%%%%%%%%%%%%%%%%%%%%%%%%%%%%%%%%%%%%%%%%%%%%%%%%%%%%%%%%%%%%%%%%%%%%
\section{Projective, injective and flat objects in $n$-trivial extension of abelian categories}
In this section we investigate the classical notions of projectivity, injectivity and flatness in the right and left $n$-trivial extension of an abelian category.\\

Throughout this part, we use the notations and conventions in Section 2, in particular whenever the category $\ntex$ is considered, it will be assumed that
$\Fi$ is right exact and accepts a right adjoint functor $\Gi$ for every $i\in I$.\\

We know that when a category $\A$ has a projective generator $P$, the class of projective objects is nothing but the class $Add(P)$ which is the
class of objects isomorphic to direct summands of direct sums of copies of $P$ in $\A$. The following theorem describes this class over $\ntex$.

%%%%%%%%%%%%%%%%%%%%%%%%%%%%%%%%%%%%%%%%%%%%%%%%%%
\begin{thm}\label{thm.Add}
Let $P$ be an object of $\A$.\\   An object $(X,f)$ of $\ntex$ is in $Add(T(P))$ if and only if $C(X,f)\in Add(P)$ and $(X,f)\cong T(C(X,f))$.
\end{thm}
\pr
The ``if" part is trivial. For simplicity $T(P):=(K,\kappa)$ and $C(X,f):=Cf$, one can notice that $\Coker((K,\kappa))=CT(P)\cong P$. Let $(X,f)$ be in $Add(T(P))$, we have the following commutative diagram
\begin{displaymath}
\xymatrix {   F_i K^{(I)}\ar[r]^{\kappa_i}\ar[d]&K^{(I)}\ar[d]\ar[r]^{\lambda}&P^{(I)}\ar[d]_{\epsilon}\ar[r]&0\\
F_i X\ar[r]_{f_i}&X\ar[r]_{\mu}&Cf\ar[r]&0}
\end{displaymath}
where the middle vertical morphism is induced from the projection $p:T(P)^{(I)}\longrightarrow (X,f)$. Thus it is clear that $\epsilon: P^{(I)}\longrightarrow Cf$ is a split epimorphism.\\
 Let $\pi:(X,f)\longrightarrow T(Cf)$ be defined by the following commutative diagram
\begin{displaymath}
\xymatrix {  (X,f)\ar[r]^{\pi}\ar[rd]^{i}&T(Cf)\\
 &T(P)^{(I)}\ar[u]_{T(\bar{\epsilon})}}
\end{displaymath}
with $i$  the canonical injection. Then, using $\bar{\epsilon}$ the retraction of $\epsilon$, we have $p\circ T(\bar{\epsilon})\circ\pi=id_{(X,f)}$.
Therefore $\pi$ is an isomorphism.\cqfd

%%%%%%%%%%%%%%%%%%%%%%%%%%%%%%%%%%%%%%%
The following description of generators comes directly.

\begin{prop}\label{Prop.gen}
An object $P$ in $\A$ is a generator if and only if $T(P)$ is a generator in $\ntex$.
 \end{prop}

Thus using the results above we deduce the following characterization of projective object in  $\ntex$ which generalizes \cite[Corollary 1.6]{TEx}.

 \begin{cor}\label{cor.proj}
 If $\A$ has enough projective objects then, the following assertions are equivalent for an object $(X,f)$ of $\ntex$:
\begin{enumerate}
\item  $(X,f)$ is projective.
\item  There exists a projective object $P$ in $\A$ such that $(X,f)$ is a direct summand of $T(P).$
\item   $C(X,f)$ is projective and $TC(X,f)\cong(X,f).$
\end{enumerate}
\end{cor}
%%%%%%%%%%%%%%%%%%%%%%%%%%%%%%%%%%%%%%

Dually we state without proof the injectivity context.

\begin{thm}\label{thm.Prod}
An object $(X,f)$ of $\conext$ is in $Prod(H(E))$ if and only if $K(X,f)\in Prod(E)$ and $(X,f)\cong HK(X,f)$.
\end{thm}

\begin{prop}\label{Prop.cogen}
An object $E$ in $\A$ is a cogenerator if and only if $H(E)$ is a cogenerator in $\conext$.
\end{prop}

\begin{cor}\label{cor.inj}
If $\A$ has enough injective objects, then the following assertions are equivalent for an object $(X,f)$ in $\conext$:\begin{enumerate}
    \item [i)]$(X,f)$ is injective.
    \item[ii)] There exists an injective object $E$ in $\A$ such that $(X,g)$ is a direct summand of $H(E)$.
\item [iii)] $K((X,f))$ is injective and $H(K(X,f))\cong (X,f)$.
\end{enumerate}
\end{cor}

%%%%%%%%%%%%%%%%%%%%%%%%%%%%%%%%%%%%%%%%%%%%%%%%%
%%%%%%%%%%%%%%%%%%%%%%%%%%%%%%%%%%%%%%%%%%%%%%%%%%%%%
Now that we have a characterization of projective and injective objects, we turn to the study of flat objects.
In \cite{TEx} the flat modules in the category of $n$-trivial extension of a ring were characterized using the Pontryagin duality between injective
and flat modules (see \cite[Proposition 1.14]{TEx}). Here we follow different approach by adopting the fact that a flat object can be seen as a direct
limits of projective objects. So we start with a general result. For this reason we assume    $\A$ to be a Grothendieck category.\\

Let us start with the following lemma.

\begin{lem}\label{lem.limit}
An object $(X,f)$ in $\ntex$ is a direct limit if $Cf:=C(X,f)$ is also a direct limit in $\A$.
\end{lem}
\pr First notice that $X\longrightarrow Cf\longrightarrow 0$ is exact. If   $Cf$ admits a direct limit, then there is a cocomplete  family
$\{d_i:D^i\longrightarrow Cf\}$ with $\{D^i\}_{i\in K}$ a family of projective objects in $\A$ and maps $(ij):D^{j}\longrightarrow D^i$  when $i<j$,
such that, for every $i\in K$, $\displaystyle{\lim_{\longrightarrow}}D^i=Cf $.
 We consider the following pullback diagram for every $i\in K$:
$$
\xymatrix { P^i \ar@{.>}[r]^{p_i}\ar@{.>}[d]_{\alpha_i} & X\ar[d]^{\qquad (*)} \\
              D^i\ar[r]^{d_i}& Cf  }
$$
We prove that   $\{p_i:{P^i}\longrightarrow X\}_{i\in K}$ is a cocomplete family in $\A$.\\
For $j\in K$, we have two maps
 $\alpha_j:P^j\longrightarrow D^j$ and $(ij):D^j\longrightarrow D^i$. Thus the diagram
\begin{displaymath}
\xymatrix { P^j \ar@{.>}[r]\ar@{.>}[d]_{(ij)\circ \alpha_i} & X\ar[d]\\
              D^i\ar[r]_{d_i}& Cf         }
\end{displaymath} commutes. Since the diagram $(*)$ is a pullback, we have a morphism $\widehat{(ij)}:P^j\longrightarrow P^i$, for every $i<j$, such that the following diagram commutes
\begin{displaymath}
\xymatrix { P^i \ar[r]\ar[d]_{\widehat{(ij)}} & X\\
              P^j\ar[ur]        }
\end{displaymath}
  We easily check that $\displaystyle{\lim_{\longrightarrow}}P^i=X$.\\
 Finally,  the family $(f_i:F_iX\longrightarrow X)_i=(\displaystyle{\lim_{\longrightarrow}}F_iP^j\longrightarrow \displaystyle{\lim_{\longrightarrow}}P^j)_i$  induces for every $i\in I$ a family $(p^j_i:F_iP^j\longrightarrow P^j)_{j\in K}$ such that   $(P^j,p^j_i)_{i\in I}$ is an object of $\ntex$.\cqfd

\par In  \cite{Cl} Beligiannis stated that a cleft extension of a Grothendieck category is also Grothendieck  if the functor $F$ preserve colimits. So in what follows $\Fi$ is supposed preserving colimits for every $i\in I$. \\

We denote by $Lim(P)$   the class of objects which are isomorphic to $\displaystyle{\lim_{\longrightarrow}}P^{(I)}$  for an object $P$ in $\A$.

\begin{thm}\label{thm.Lim} For an object $P$ in $\A$, an object $(X,f)$ is in $ Lim(T(P))$ if and only if $C(X,f)$  is in $Lim(P)$ and $T(C(X,f))\cong (X,f)$.
\end{thm}
\pr
Let $(X,f)=\displaystyle{\lim_{\longrightarrow}}(T(P^{(I)}))$. By the construction of $T$  we get $(X,f)=T\displaystyle{\lim_{\longrightarrow}}P^{(I)}$. Then, using the fact that $CT=id_{\A}$, we get the result. The ``if" part is easy to check.\cqfd

We know that an object in a  category is flat if it is a direct limit of projective objects. Then, a characterization of flat object follows directly from the  results above.

\begin{cor} \label{cor.flat}
An object $(X,f)$ in $\ntex$ is flat if and only if  $C((X,f))$ is flat and $T(C(X,f))\cong (X,f)$.
\end{cor}

As a consequence we establish a characterization of flat modules over $\Rn$ which generalizes \cite[Proposition 1.14]{TEx}.

\begin{cor}\label{cor.flatmod}
A left  $\Rn$-module   $(X,f)$ is flat if and only if $\Coker (X,f)$ is flat and the following sequence $(**)$
\begin{displaymath}
 \xymatrixcolsep{3pc}\xymatrix{\oplus_{j}\oplus_{i}(M_j\otimes M_i\otimes X)\ar[r]^-{h}& \oplus_{i}(M_iX)\ar[r]^-{(f_n,..., f_1)} & X\ar[r]& \Coker(X,f)\ar[r]&0 .  }
\end{displaymath}
is exact  for every $i,j\in J$, where  $$h=\left(
 \begin{array}{cccc}
   M_n\otimes f_n,... ,M_n\otimes f_1 &\ldots &0\\
                      0 &\ddots &0\\
                 \vdots &\ddots &\vdots\\
                        0&\cdots & M_1\otimes f_n,... ,M_1\otimes f_1\\
                       \end{array}
                 \right).   $$
\end{cor}
\pr Just remark that the exactness of the sequence $(**)$ is equivalent to the fact that $T(C(X,f))\cong (X,f)$.\cqfd

%%%%%%%%%%%%%%%%%%%%%%%%%%%%%%%%%%%%%%%%%%%%%%%%%%%%%%%%%%%%%
%%%%%%%%%%%%%%%%%%%%%%%%%%%%%%%%%%%%%%%%%%%%%%%%%%%%%%%%%%%%%%%%%%%%%
%%%%%%%%%%%%%%%%%%%%%%%%%%%%%%%%%%%%%%%%%%%%%%%%%%%%%%%%%%%%%%%%
%%%%%%%%%%%%%%%%%%%%%%%%%%%%%%%%%%%%%%%%%%%%%%%%%%%%%%%%%%%%%%%%%%
\section{Applications}
This section is devoted  to the study of some properties of $\Rn$ as  applications of the study done in the previous section. \\

We first show when $\Rn$ is $m$-perfect. A ring $R$ is said to be  $m$-perfect for some positive integer $m$ if every flat module has projective dimension at
most $m$ (see  \cite[Definition 1.1]{n-perfect}). We need the following lemma.

\begin{lem}\label{lem.pdT}
Let $(X,f)$ be an object in $\ntex$ such that  $X=X_1\oplus X_2$ with $\Im(f_i)\in X_2$, for all $i\in I$, then
$$\pda(X_1)\leq \pde(X,f). $$
especially we have
 $ \pda X \leq \min(\pde Z(X), \pde T(X)) $  for every object $X\in \A$.
\end{lem}
\pr
We assume that $\pde(X,f)=m<\infty$. For $m=0$, we have $(X,f)$ is a projective object of $\ntex$. Hence $C(X,f)$ is also projective in $\A$ by Corollary  \ref{cor.proj}. Since $\Im(f_i)\in X_2$ and $X=X_1\oplus X_2$ we get that $X_1$ is a direct summand of $C(X,f)$. Therefore $X_1$ is projective in $\A$.\\
 Now assume that $m\geqslant 1$. Consider two exact  sequences in $\A$ $$\xymatrix{
       0\ar[r]&K\ar[r] &P_1\ar[r]^{\epsilon_1}&X_1\ar[r]&0\\                          } \quad \mathrm{and}\quad \xymatrix{
       P_2\ar[r]^{\epsilon_2}&X_2\ar[r]&0}$$  with $P_1$ and $P_2$ projective objects in \A. Since $\Im f_i\in X_2$ we can set  $f_i=(f_{i,1},f_{i,2})$ with  $f_{i,1}: \Fi X_1\longrightarrow X_2$ and
$f_{i,2}=\Fi X_2\longrightarrow X_2$. Now, let $L$ be the kernel of the morphism:
$$\lambda=((f_{i,1}\circ\Fi\epsilon_{1})_{i\in I},\epsilon_2,(f_{i,2}\circ\Fi\epsilon_2)_{i\in I}):\displaystyle\bigoplus_{i\in I}\Fi P_1\oplus P_2\displaystyle\bigoplus_{i\in I}\Fi P_2\longrightarrow X_2$$
and consider
$\mu=\left( \begin{array}{cc}
 \epsilon_1 &0\\
 0 & \lambda
\end{array}\right):P_1  \displaystyle\bigoplus_{i\in I}\Fi P_1\oplus P_2\displaystyle\bigoplus_{i\in I}\Fi P_2\longrightarrow X_1\oplus X_2.$\\
Then we obtain the following exact sequence:
\begin{displaymath}
    \xymatrix{
       (K\oplus L,\alpha)\ar[r]^-{\Ker\mu} &T(P_1)\oplus T(P_2)\ar[r]^-{\mu}&(X,f_i)\ar[r]&0\\ }
\end{displaymath}
with $\alpha=(\alpha_{i,1}: \Fi K\longrightarrow L,\alpha_{i,2}:\Fi L\longrightarrow L )$.
The middle term represents a projective object in $\ntex$, hence $\pde(X,f_i)=1+\pde(K\oplus L,\alpha_i) $. By the induction  $\pde(K\oplus L,\alpha_i)\geq \pda K $. Therefore $\pde(X,f)\geq 1+\pda K=\pda X_1.$ \cqfd

Now we can give the first result which generalizes \cite[Proposition 1.15]{TEx}.

\begin{thm}\label{thm.perfect}
The $n$-trivial extension $\Rn$ is $k$-perfect (for $k\in \mathbb{N}$) if and only if $R$ is $k$-perfect.
\end{thm}
\pr For simplicity, the Pontryagin dual $\Hom_Z(X, \mathbb{Q}/ \mathbb{Z}) $ of a module $X$ will be denoted by $X^*$.
  We suppose that $R$ is $k$-perfect. Let $(X,f)$ represent a flat object. By Corollary \ref{cor.flatmod}, $C=\Coker(X,f)$ is flat, hence there exists an exact projective resolution of $R$-modules, \begin{displaymath}
 \xymatrix{0\ar[r]& P_k\ar[r]& P_{k-1}\ar[r]& \cdots \ar[r]&P_1\ar[r]&C\ar[r]&0.}
\end{displaymath} By Corollary \ref{cor.flatmod}, we  have  $X^{*}\cong C^{*}\displaystyle \bigoplus_{i\in J} \Gi C^{*}  $ and $ \Gi C^{*}\cong (\Fi C)^* $, then  $(X,f)\cong T(C)$.\\
Since $C$ is flat, we get an exact sequence, \begin{displaymath}
 \xymatrix{0\ar[r]& T(P_k)\ar[r]& T(P_{k-1})\ar[r]& \cdots \ar[r]&T(P_1)\ar[r]&T(C)\ar[r]&0.}
 \end{displaymath}
Therefore, $\Rn$  is $k$-perfect.\\
Conversely, let $X$ be a flat $R$-module, by \ref{cor.flatmod} and the fact that $CT(X)\cong X$, $T(X)$ is also flat. By Lemma \ref{lem.pdT},
we have $\pda X \leq \pde T(X)\leq k$ since  $\Rn$ is $k$-perfect. Therefore, $R$ is
$k$-perfect.\cqfd

We end the paper with a result that computes the injective dimension of the ring
$\Rn$. For this we need the following lemma.

\begin{lem}\label{lem.inj}
If $X$ is an object in $\A$ such that $R_i\Gj X=0$ for all $ j\in I$ and $ i>0$, then $$\ide(HX)=\ida X. $$
\end{lem}
\pr
Consider an injective resolution of $X$ in $\A$
\begin{displaymath}
    \xymatrix{
       0\ar[r]&X\ar[r] &I_1\ar[r]&I_2\ar[r]&\cdots.\\                          }
\end{displaymath}
The hypothesis $R_i\Gj X=0$ shows that the following sequence is an injective resolution  of $H(X)$ in $\conext$ \begin{displaymath}
    \xymatrix{
       0\ar[r]&H(X)\ar[r] &H(I_1)\ar[r]&H(I_2)\ar[r]&\cdots.                         }
\end{displaymath}
Hence, $\ida(X)\geq\ide H(X)$.\\
The dual proof of Lemma \ref{lem.pdT} gives $\ida(X)\leq\ide H(X)$. \cqfd

The following result generalizes \cite[Theorem 4.32]{TEx}.

\begin{thm}\label{thm.selfinj}
 \par Suppose that:
\begin{equation*} (*)\quad \Hom(M_i,M_n)=
                    \begin{cases}
                 R & \text{ for } i=n.\\
              M_{n-i} &\text{ if } i<n.\\
                           \end{cases}
                      \end{equation*}
and $\Ext^{k}(M_i,M_n)=0$ for $k\geq 1.$
Then $ \idrn(R\ltimes_{n} M)=\idr (M_n).$
\end{thm}
\pr From Lemma \ref{lem.inj}, $H(M_n)$ has injective dimension equal to $\idr M_n$.\\
 We remark that $\Rn$ is isomorphic to the object $T(R)=(R\oplus(\displaystyle\bigoplus_{i\in I}F_i R), \kappa)$. Under the category isomorphism $\lRM \cong\RRM$, the object $T(R)$ corresponds to the object
  $\widetilde{T(R)}=(R\oplus(\displaystyle\bigoplus_{i\in I}F_i R), \widetilde{\kappa})$.
   By the following equalities   $G_iM_n=M_{n-i}$,
   $F_iR=M_i$ and $G_n M_n=R$,   $\widetilde{T(R)}$ is isomorphic to $H(M_n)$ and hence $ \idrn(R\ltimes_{n} M)=\idr(M_n).$\cqfd

%%%%%%%%%%%%%%%%%%%%%%%%%%%%%%%%%%%%%%%%%%%%%%%%%%%%%%%%%%%
%%%%%%%%%%%%%%%%%%%%%%%%%%%%%%%%%%%%%%%%%%%%%%%%%%%%%%%%%

\noindent {\bf Acknowledgement.}  
The authors would like to thank the referee for the careful reading.\\
Dirar Benkhadra's research reported in this publication was supported by a scholarship from the Graduate Research Assistantships in Developing Countries Program of the Commission for Developing Countries of the International Mathematical Union.

%%%%%%%%%%%%%%%%%%%%%%%%%%%%%%%%%%%%%%%%%%%%%%%%%%%%%%%%%%%
%%%%%%%%%%%%%%%%%%%%%%%%%%%%%%%%%%%%%%%%%%%%%%%%%%%%%%%%%
%%%%%%%%%%%%%%%%%%%%%%%%%%%%%%%%%%%%%%%%%%%%%%%%%%%%%%%%%
%%%REFERENCES%%%%%%%%%%%%%%%%%%%%%%%%%%%%%%%%%%%%%%%%%%%%
%%%%%%%%%%%%%%%%%%%%%%%%%%%%%%%%%%%%%%%%%%%%%%%%%%%%%%%%

\bigskip\bigskip

%%%%%%%%%%%%%%%%%%%%%%%%%%%%%%%%%%%%%%%%%%%%%%%%%%%%%%%%

\end{document}